\documentclass[12pt]{amsart}
\usepackage{amsfonts,amssymb,amscd,amsthm,amsmath,euscript}
\usepackage{graphpap}

 
\newtheorem{theorem}[subsection]{Theorem}
\newtheorem*{theorem*}{Theorem}
\newtheorem{lemma}[subsection]{Lemma}
\newtheorem{proposition}[subsection]{Proposition}
\newtheorem{corollary}[subsection]{Corollary}

\theoremstyle{definition}
\newtheorem{definition}[subsection]{Definition}

\newtheorem{example}[subsection]{Example}
\theoremstyle{remark}
\newtheorem{remark}[subsection]{Remark}

\theoremstyle{definition}
\newtheorem{subs}[subsection]{}

\newcommand{\mt}[1]{\operatorname{#1}}

\newcommand{\DDD}{{\mathbb D}}
\newcommand{\AAA}{{\mathbb A}}
\newcommand{\QQ}{{\mathbb Q}}
\newcommand{\ZZ}{{\mathbb Z}}
\newcommand{\CC}{{\mathbb C}}
\newcommand{\OO}{{\mathcal O}}

\newcommand{\PP}{{\mathbb P}}

\newcommand{\NN}{{\mathbb N}}
\newcommand{\FFF}{{\mathbb F}}
\newcommand{\PPP}{{\EuScript{P}}}
\newcommand{\an}{{\cong_{\mt {an}}}}
\newcommand{\SM}{{\Phi_{\bf sm}}}
\newcommand{\M}{{\Phi_{\bf m}}}
\newcommand{\Supp}{\mt{Supp}}
\newcommand{\Pic}{\mt{Pic}}
\newcommand{\Sing}{\mt{Sing}}
\newcommand{\Diff}{\mt{Diff}}

\newcommand{\Exc}{\mt{Exc}}

\newcommand{\num}{\nu_{\mt{num}}}
\newcommand{\NE}{\overline{\mt{NE}}}
\newcommand{\down}[1]{\llcorner #1 \lrcorner}
\newcommand{\downn}[1]{\left[ #1\right]}

\newcommand{\fr}[1]{\{ #1\}}

\title{Complements on log surfaces}
\author{S.~A.~Kudryavtsev}

\date{}
\address{Department of Algebra, Faculty of Mathematics,
Moscow State Lomonosov University, 117234 Moscow, Russia}

\email{kudryav@mech.math.msu.su\newline\
Web page: www.math.msu.su/department/algebra/staff/kudryav.htm}

\begin{document}
\begin{abstract} More strong version of the
main inductive theorem about the complements on surfaces is proved and
the models of exceptional log del Pezzo surfaces with
$\delta=0$ are constructed\thanks{This work was done with the partial support of the
Russian Foundation for Basic Research (grant no. 02-01-00441), the Leading
Scientific Schools (grant no. 00-15-96085) and INTAS-OPEN
(grant no. 2000\#269.)}.
\end{abstract}
\maketitle

\section*{\bf {Introduction}}
The theory of complements on algebraic varieties has been created by
V.~V.~Shokurov in the papers \cite{Sh1}, \cite{Sh2}.
It is a powerful tool for studying algebraic varieties, extremal contractions
and singularities.
Roughly speaking, the complement is a "good"\ divisor in the multiple anticanonical
linear system.
The advantage of this theory is that the concept of complement is an invariant
in Log Minimal Model Program.
Moreover a complement has an inductive property, this means that
the complement
finding problem for an
$n$-dimensional variety is reduced to the same one for an
$(n-1)$-dimensional variety.
See the papers \cite{Sh2}, \cite{PrLect},
\cite{PrSh} with reference to the theory of complements on the
high-dimensional varieties.
For example, the application of this theory for the three-dimensional varieties
is given in the
papers \cite{Kud3}, \cite{Kud4}, \cite{PrRed}, \cite{PrBound}.
\par
Thus, in order to study effectively
the three-dimensional contractions and singularities
it is important to classify the log del Pezzo surfaces completely.
The last open two-dimensional problem (in the framework of the theory of complements)
is the classification of exceptional log del Pezzo surfaces.
The exceptional log del Pezzo surfaces $(S,D)$ are divided into three types:
$\delta(S,D)=0$, $\delta(S,D)=1$ and
$\delta(S,D)=2$, where
\begin{eqnarray*}
\delta(S,D)=\#\Big\{E\mid E\ \text{is a divisor with a discrepancy}\
a(E,D)\le -\frac67   \Big\}.
\end{eqnarray*}
The cases $\delta(S,D)=1$ and $\delta(S,D)=2$ were classified in the papers
\cite{KudLd}, \cite{Sh2}.
To study the remaining case
$\delta(S,D)=0$
the theory of complements on surfaces must be applied in more wide set of
coefficients.
Therefore it will be considered when all coefficients of a boundary
$D$ are greater then or equal to 1/2.
\par
One of the main results of this paper given in
\S 2 is more strong version of the
main inductive theorem about the complements on surfaces.
Using this inductive theorem we construct the models of
exceptional log del Pezzo surfaces with $\delta=0$ in
\S 4 (see
definition \ref{maindef}). In \S 3 we give the classification of non-rational
exceptional log surfaces.
Also, one type of exceptional log del Pezzo surfaces
with $\delta=0$ is described completely in \S 4 (see
theorem \ref{main3}).
\par
I am grateful to Professor  Yu.G.~Prokhorov for valuable remarks.

\section{\bf {Preliminary facts and results}}
All varieties are algebraic and are assumed to be defined over
$\CC$, the complex number field.
The main definitions, notations and notions used in the paper are
given in \cite{Koetal}, \cite{PrLect}.

\begin{definition}
Put $\SM=\{1-1/m\mid m\in\ZZ_{>0}\cup\{\infty\}\}$ and
$\M=\SM\cup [6/7,1]$. A coefficient $d$ is called {\it
standard} if $d\in\SM$.
\par
Put $\ZZ/(n)=\{k/n\mid k\in\ZZ_{>0}\}$.
\end{definition}

\begin{definition} For fix $n\in \NN$ put
\begin{eqnarray*}
\PPP_n=\{ a \mid 0\le a \le 1,\ \down{(n+1)a}\ge n a \}.
\end{eqnarray*}
It is clear that
\begin{multline*}
\PPP_n=\left\{0\right\}\bigcup \downn{\frac1{n+1},\frac1{n}}
\bigcup \downn{\frac2{n+1},\frac2{n}}\bigcup\dots\\
 \bigcup\downn{\frac k{n+1},\frac k{n}}\bigcup\dots\bigcup
\downn{\frac{n}{n+1},1}.
\end{multline*}
\end{definition}

\begin{definition}
Let us define the set $\Phi^{[a,b)}=\big(\SM\cap [0,b)\big)\cup
[a,b).$ Put $\Phi_i=\Phi^{[a,b)}$, where
$a=\frac{i-1}{i}$ and $b=\frac{i}{i+1}$.
\end{definition}

\begin{definition}
The pair $(X,D=\sum d_kD_k)$ is called {\it a pair of type} $\Phi_i$ if
the three following conditions are satisfied:
\begin{enumerate}
\item $d_k\in \Phi_i$ for all $k$;
\item there exists $j$ such that
$d_j\ge \frac{i-1}{i}$;
\item the divisor $K_X+D$ is
$\frac1i$-log terminal.
\end{enumerate}
\end{definition}

\begin{definition}
Let $(X/Z\ni P,D)$ be a pair, where $D$ is a subboundary. Then a
\textit{$\QQ$-complement} of $K_X+D$ is a
log divisor $K_X+D'$ such that $D'\ge D$, $K_X+D'$ is log canonical and
$n(K_X+D')\sim 0$ for some $n\in\NN$.
\end{definition}

\begin{definition}\label{def1}
Let $X$ be a normal variety and let $D=S+B$ be a subboundary on
$X$ such that $B$ and $S$ have no common components, $S$ is an
effective integral divisor and $\down{B}\le 0$. Then we say that the divisor
$K_X+D$ is \textit{$n$-complementary} if there is a $\QQ$-divisor
$D^+$ such that
\begin{enumerate}
\item
$n(K_X+D^+)\sim 0$ (in particular, $nD^+$ is an integral divisor);
\item the divisor
$K_X+D^+$ is log canonical;
\item
$nD^+\ge nS+\down{(n+1)B}$.
\end{enumerate}
In this situation the
\textit{$n$-complement} of
$K_X+D$ is $K_X+D^+$. The divisor $D^+$ is called an
\textit{$n$-complement} too.
\par
Let $X$ be a semi-smooth variety in codimension 1. Then the divisor
$K_{X}+D$ is $n$-{\it semicomplementary} if there is a
$\QQ$-divisor $D^{+}$ satisfying previous conditions
(1), (3) and next condition $(2')$:
\begin{enumerate}
\item[$(2')$] the divisor $K_{X}+D^+$ is semi log canonical.
\end{enumerate}
\end{definition}

\begin{proposition}\cite[example 5.2]{Sh1}, \cite[theorem 19.4]{Koetal}
\label{compl1} Let $X$ be a semi-smooth, connected, complete curve.
Let $D$ be a boundary on $X$ contained in the smooth part of $X$.
Assume that
$-(K_X+D)$ is a nef divisor on every component of $X$.
Then
\begin{enumerate}
\item the divisor $K_X+D$ is $n$-semicomplementary for
$n\in\{1,2,3,4,6\}$;
\item if the divisor $K_X+D$ is not $1$- and
$2$-semicomplementary then $X\cong\PP^1$ and $\down{D}=\down{D^+}=0$.
\end{enumerate}
\end{proposition}

The following statements show the invariance of complements with respect to
the log minimal model program and their inductive properties.
\begin{proposition}\cite[lemma 5.4]{Sh1}
Let $f\colon X\to Z$ be a birational contraction of varieties and let
$D$ be a subboundary on $X$. If the divisor $K_X+D$ is
$n$-complementary then the divisor $K_Z+f(D)$ is also $n$-complementary.
\end{proposition}

\begin{proposition}\label{compl3}
\cite[lemma 4.4]{Sh2} Let $f\colon X\to Z$ be a birational contraction of
varieties and let $D$ be a subboundary on $X$.
Assume that
\begin{enumerate}
\item the divisor $K_X+D$ is $f$-nef;
\item the coefficient of every non-exceptional component of
$D$ meeting
$\Exc f$ belongs to $\PPP_n$;
\item the divisor $K_Z+f(D)$ is
$n$-complementary.
\end{enumerate}
Then the divisor $K_X+D$ is also $n$-complementary.
\end{proposition}

For the two-dimensional varieties we have more strong theorem about the
inductive property of complements
then for the high-dimensional varieties
\cite[proposition 4.4.1]{PrLect}.
\begin{theorem}\label{compl2} \cite[proposition 4.4.3]{PrLect}
Let $(X/Z\ni P,D=S+B)$ be a log surface with the following properties:
\begin{enumerate}
\item the divisor $K_X+D$ is divisorial log terminal;
\item the divisor $-(K_X+D)$ is nef and big over $Z$;
\item $S=\down{D}\ne 0$ in the neighborhood of $f^{-1}(P)$.
\end{enumerate}
Assume that near $f^{-1}(P)\cap S$ there exists an
$n$-semicomplement $K_S+\Diff_S(B)^+$ of $K_S+\Diff_S(B)$. Then near
$f^{-1}(P)$ there exists an $n$-complement $K_X+S+B^+$
of $K_X+S+B$ such that $\Diff_S(B)^+=\Diff_S(B^+)$.
\end{theorem}

\begin{definition}
Let $(X/Z\ni P,D)$ be a contraction of varieties, where $D$ is a boundary.
In the case when
$\dim Z\ne 0$ the contraction is said to be {\it exceptional}
if for every
$\QQ$-complement $D'$ there is at most one divisor $E$ (not
necessarily exceptional) such that $a(E,D')=-1$. In the case when
$\dim Z=0$ the log variety is said to be {\it exceptional} if the pair
$(X,D')$ is kawamata log terminal for every
$\QQ$-complement $D'$.
\end{definition}

\begin{definition}
Let $(X,D)$ be an exceptional log variety. Define
\begin{eqnarray*}
\delta(X,D)=\#\Big\{E\mid E\ \text{is an exceptional or non-exceptional divisor}\\
\text{with a discrepancy}\
a(E,D)\le -\frac67   \Big\}.
\end{eqnarray*}
\end{definition}

\begin{lemma}\label{two1}
Let $(X\ni P, \alpha C+B)$ be a germ of two-dimensional log terminal pair, where
 $(X\ni P)$ is a non-cyclic singularity,
$C$ is a curve, $B\ge 0$ and $\alpha\ge 0$. Then
\begin{enumerate}
\item the divisor $K_X+\alpha C+B$ is not
$(1-\alpha)$-log terminal;
\item the divisor $K_X+\alpha C+B$ is strictly
$(1-\alpha)$-log canonical if and only if $(X\ni
P,\alpha C+B)\an (\CC^2\ni 0,\alpha\{xy=0\})/\DDD_n$, where
$\DDD_n\subset SL_2(\CC)$ is a dihedral subgroup.
\end{enumerate}
\begin{proof} For some number $0<c\le 1$ the pair $(X,cC+B)$ is log canonical,
but not purely log terminal
\cite[theorem 2.1.2]{PrLect}. Let
$f\colon (Y,E)\to (X\ni P)$ be an inductive blow-up of this log pair
(\cite[theorem 1.9]{Kud2}, \cite[proposition 3.1.4]{PrLect}).
Then the divisor $K_Y$ is $f$-nef since $f$ is the blow-up of the central
vertex of  minimal resolution graph
\cite[\S 6]{PrLect}. Therefore we have
$$
a(E,\alpha C+B)=-\frac{\alpha}c+(1-\frac{\alpha}c)\cdot a(E,B)\le
-\frac{\alpha}c \le -\alpha.
$$

The equality holds if and only if
$c=1$, $B=0$ and $a(E,0)=0$. By the classification of two-dimensional log terminal
singularities we obtain the required statement
(for example, see \cite[theorem 2.1.2]{PrLect}).
\end{proof}
\end{lemma}

\begin{lemma}\label{two2}
Let $(X\ni P, \alpha C+B)$ be a germ of two-dimensional log terminal pair, where
$(X\ni P)$ is a cyclic singularity,
$C$ is a curve, $B=\sum b_iB_i\ge 0$ and $\alpha\ge 0$. Assume that
the pair $(X,C+B)$ is not purely log terminal and $b_i\ge \frac12$ for all
$i$. Then
\begin{enumerate}
\item the divisor $K_X+\alpha C+B$ is not
$(1-\alpha)$-log terminal;
\item the divisor $K_X+\alpha C+B$
is strictly $(1-\alpha)$-log canonical if and only if:
\begin{enumerate}
\item $(X\ni P,\alpha C+B)\an (\CC^2\ni
0,\alpha\{xy=0\})/\ZZ_n(n-1,1)$;
\item $(X\ni P,\alpha C+B)\an
(\CC^2\ni 0,\alpha\{x^2+y^2=0\})/\ZZ_4(3,1)$;
\item $(X\ni
P,\alpha C+B)\an (\CC^2\ni 0,\alpha\{x^2+y^4=0\})/\ZZ_2(1,1)$, where
$\alpha\le 2/3$;
\item $(X\ni P,\alpha C+B)\an (\CC^2\ni
0,\frac12\{x=0\}+\frac12\{x+y^3=0\})/\ZZ_2(1,1)$.
\end{enumerate}
\end{enumerate}
\begin{proof}
Assume that $C$ is a reducible curve or $B$ consists of
at least two divisors. Let $\psi\colon\widetilde X\to X$
be a minimal resolution and $\Gamma$ be its graph. The proper transforms of
$C$ and $B$ are denoted by $\widetilde C$ and $\widetilde B$. The curve of $\Gamma$
intersecting $\widetilde C$ is denoted by $\widetilde E$.
Let us contract all other curves of
$\Gamma$. We obtain a blow-up
$f\colon (Y,E)\to (X\ni P)$. Write
$$
K_Y+aE+B_Y+C_Y=f^*(K_X+C+B).
$$
Since
\begin{equation}
0\le -2+\deg\Diff_E(0)+(B_Y+C_Y)\cdot E=(1-a)E^2
\end{equation}
then $a\ge 1$. Hence
\begin{equation}
a(E,\alpha C+B)=-\alpha a+(1-\alpha)a(E,B)\le -\alpha a\le
-\alpha.
\end{equation}

The equality holds if and only if
$a=1$, $B=0$ and $a(E,0)=0$.
By the classification of two-dimensional log terminal
singularities we obtain subcase
(a), or subcase (d) considered below (in this situation
$B=0$).
\par
Assume that $C$ is an irreducible curve and $B$ consists of at most one
divisor. If $(X\ni P)\ncong_{\rm an}(\CC^2\ni 0)/\ZZ_n(1,1)$ then
arguing as above we can find the curve $E$ such that inequality
(1) holds. If we have an equality in
(2) then we obtain subcase (b).
\par Let
$(X\ni P)\an(\CC^2\ni 0)/\ZZ_n(1,1)$. Put $E=\Exc \psi$. If
$(\widetilde C+\widetilde B)\cdot E\ge 2$ then lemma is proved by the same arguments.
If we have an equality in
(2) then we obtain subcase
(c). Therefore we may assume that
$$(X\ni P,C)\an (\CC^2\ni 0,\{x=0\})/\ZZ_n(1,1),
B=b_1B_1\ \text{and}\ \widetilde C\cap \widetilde B_1\cap E=P.$$

Take an usual blow-up at the point $P$. Then
$$
a(E',\alpha
C+B)=-\alpha\big(1+\frac1n\big)-b_1\big(1+\frac1n\big)+\frac2n,
$$
where $E'$ is a corresponding exceptional curve.
It can easily be checked that
$a(E',\alpha C+B)>-\alpha$ if and only if $n=2$, $\alpha<1/2$ and we have
$a(E,\alpha C+B)\le-\alpha/2-1/4< -\alpha$. Moreover
$a(E',\alpha C+B)=-\alpha$ if and only if $n=2$,
$\alpha=b_1=\frac12$, that is, we obtain subcase (d) (in this situation $B\ne 0$).
\end{proof}
\end{lemma}

\begin{proposition}\label{two3}
Let the pair $(X\ni P,\alpha C+B)$ be of type $\Phi_i$, where
$i=2,3,4,5,6$ and $\alpha\ge \frac{i-1}i$. Then one of the following possibilities
holds.
\begin{enumerate}
\item[1)] $(X\ni P,\alpha C+B)\an(\CC^2\ni
0,\alpha\{x=0\}+\beta\{x+y^k=0\})$, where $i=2$, $k\ge 2$,
$\alpha>1/2$, $\beta>1/2$, $\alpha+\beta<1+\frac1{2k}$.
\item[2)]
$(X\ni P,\alpha C+B)\an(\CC^2\ni
0,\alpha\{x=0\}+\frac12\{x+y^k=0\})$, where $k\ge 2$. If $k=2$ then
$i$ is arbitrary. If $k=3$ then $i=2,3$. If $k\ge 4$ then $i=2$ and
$\alpha<\frac12+\frac1{2k}$.
\item[3)] $(X\ni P,\alpha
C+B)\an(\CC^2\ni 0,\alpha\{x^2+y^k=0\})$, where $k\ge 2$. If
$k=2$ then $i$ is arbitrary. If $k=3$ then $i=2,3,4$. If $k\ge 4$ then
$i=2$ and $\alpha<\frac12+\frac1{4\down{k/2}}$.
\item[4)] $(X\ni
P,\alpha C+B)\an(\CC^2\ni 0,\alpha\{x=0\}+b_1\{y=0\})/\ZZ_n(q,1)$,
where $(n,q)=1$ and $\frac{\frac{n}{i}-1+b_1}{1-\alpha}<q\le n$
$($the case $n=q=1$ is possible$)$.
\end{enumerate}
\begin{proof} By lemma \ref{two1} $(X\ni P)$ is a cyclic singularity or a
smooth point.
The condition that $K_X+D$ is $\frac1{i}$-log terminal divisor and the form
of the coefficients of a divisor
$D$ are principal in the proposition proof.
\par
Assume that the divisor $K_X+C+B$ is purely log terminal. If
$(X\ni P)$ is a smooth point then $(X\ni P,\alpha C+B)\an(\CC^2\ni
0,\alpha\{x=0\}+b_1\{y=0\})$.
By the same argument as in the proof of proposition
1.9 \cite{KudLd} we obtain case 4).
\par Assume that the divisor $K_X+C+B$ is not purely log terminal. Then by lemma
\ref{two2} $(X\ni P)$ is a smooth point. If the divisor
$K_X+C+B$ is divisorial log terminal then we obtain case
3) with $k=2$. Suppose that the divisor
$K_X+C+B$ is not divisorial log terminal.
Then there are two possibilities for a divisor
$B$.
\par Let $B=b_1B_1\ne 0$. Then $C$ is a smooth curve and $B_1$ is tangent $C$
at the point $P$. Therefore we obtain cases 1) and 2).
\par Let $B=0$. Then $C$ is a singular curve and we obtain case 3) with $k\ge 3$.
\end{proof}
\end{proposition}

\begin{remark} The case $\M=\bigcup^{\infty}_{i\ge 7}\Phi_i$ is developed
in proposition 1.9 \cite{KudLd}.
\end{remark}

\section{\bf {Main inductive theorem on surfaces}}
Next theorem \ref{mainind} is more strong version of the
main inductive theorem about the complements on surfaces
\cite[theorem 2.3]{Sh2}.
\begin{theorem}\label{mainind}
Let $(S,D=\sum d_iD_i)$ be a projective log surface with the following properties:
\begin{enumerate}
\item the divisor $K_S+D$ is log canonical, but not kawamata log terminal;
\item the divisor $-(K_S+D)$ is nef;
\item there exists a $\QQ$-complement of $K_S+D$;
\item $d_i\ge \frac12$ for all $i$.
\end{enumerate}
\par
Then there is $1$-, $2$-, $3$-, $4$- or $6$-complement
of $K_S+D$ which is not kawamata log terminal,
except the cases from example $\ref{mainind1}$.
\par
Besides, if there is an infinite number of divisors
$E$ with a discrepancy
$a(E,D)=-1$ then there is $1$- or $2$-complement of
$K_S+D$ which is not kawamata log terminal.
\begin{proof}
In many cases this theorem is true without condition
(4). When proving this theorem we follow the paper
\cite{Sh2}. The cases using condition
(4) are considered in details.
\par
Applying a minimal log terminal modification
\cite[definition 3.1.3]{PrLect} we may assume that the pair $(S,D)$
is divisorial log terminal. Put
$C=\down{D}\ne 0$ and $B=\fr{D}$.
We have three cases depending on the numerical dimension of a divisor
$-(K_S+C+B)$.
\par {\bf Case I.} Assume that $-(K_S+C+B)$ is a big divisor.
Then all required statements immediately follow by
proposition \ref{compl1} and theorem \ref{compl2}.
Let us remark that condition
(4) on a boundary $D$ is unnecessary in this case.
\par Before discussing two remaining cases let us make more precise the
structure of a log surface $(S,D)$.
\par Let $S$ be a non-rational surface.
Then our theorem is proved in
\cite[theorem 2.3]{Sh2}, \cite[theorem 8.2.1]{PrLect}. Moreover
there exists
1- or 2-complement which is not kawamata log terminal
and there are at most two divisors
$E$ with a discrepancy $a(E,D)=-1$.
Let us remark that condition
(4) on a boundary $D$ is also unnecessary in this case.
\par Let $C$ be not the chain of rational curves.
Then our theorem is also true without condition
(4) on the coefficients of a boundary
$D$ \cite[theorem 2.3]{Sh2}.
\par Thus we may assume that
$S$ is a rational surface and $C$ is a chain of rational curves.
\par {\bf Case II.} Assume that $K_S+C+B\not\equiv 0$ and
$-(K_S+C+B)$ is not a big divisor. By proposition 2.5
\cite{Sh2} we can assume that the divisor $-(K_S+C+B)$ is semi-ample.
Let $\nu\colon S\to Z\cong\PP^1$ be the morphism given by a linear system
$|-m(K_S+C+B)|$, where $m\gg 0$. The next lemma is basic to construct the
complements.
\begin{lemma}\cite[lemma 2.20, lemma 2.21]{Sh2}\label{compl4}
Let $(S,D)$ be a projective log surface with a structure of fibration
onto a curve
$f\colon S\to Z$, where $D$ is a boundary. Let
$C=\down{D}\ne 0$ and $B=\fr{D}$.
Assume that the following conditions are satisfied:
\begin{enumerate}
\item there exists a section $C_1\subset C$ of $f$;
\item the divisor
$K_C+\Diff_C(B)$ is $n$-semicomplementary;
\item the divisor
$K_S+C+\down{(n+1)B}/n$ is numerically trivial on a general fiber;
\item
the divisor $-(K_S+C+B)$ is nef;
\item the divisor $K_S+C+B$
is log terminal in some $($analytic$)$ neighborhood of a divisor $C$.
\end{enumerate}
Then the divisor $K_S+C+B$ is $n$-complementary.
Moreover conditions
$(1)$ and $(3)$ can be replaced by condition $(1'):$
\begin{enumerate}
\item[$(1')$] there exists a multi-section $C_1\subset C$ of $f$ and
$S$ is a rational surface.
\end{enumerate}
\end{lemma}
There are three possibilities for $C$.
\par A). Let $C_1\subset C$ be a multi-section of $\nu$.
Then the required statements don't depend on condition
(4) on the coefficients of $D$ and follow by lemma \ref{compl4}
and proposition \ref{compl1}.
\par B). Let $C$ has the unique section $C_1$ of $\nu$.
Lemma \ref{compl4} cannot be applied
if and only if there is a horizontal component $B_i$ of
$B$ with a coefficient
$b_i\in \ZZ/(n+1)$ \cite[lemma 2.27]{Sh2}.
\par
On the other hand,
if there is a horizontal component $B_j$ of
$B$ with a coefficient $b_j\not\in\ZZ/(n+1)$,
then we consider the divisor
$K_S+C+B-\varepsilon B_j$, where $0<\varepsilon\ll 1$.
It has the same
$n$-complements as the divisor $K_S+C+B$ (see definition \ref{def1}).
Since $-(K_S+C+B)+\varepsilon B_j$ is a nef and big
divisor then our theorem is reduced to case I.
Therefore we assume that all horizontal components of
$B$ have the coefficients from the set $\ZZ/(n+1)$.
\par
Assume that $C$ is a reducible divisor.
Then the divisor $K_C+\Diff_C(B)$ is 1- or 2-semicomplementary
by proposition \ref{compl1}.
\par
If it is 1-semicomplementary then $B^{\mt{hor}}=\frac12B_1+\frac12B_2$ or
$B^{\mt{hor}}=\frac12B_1$.
Hence the divisor $K_S+C+B$ is 2-semicomplementary.
\par
If it is 2-semicomplementary then we have a contradiction with condition (4).
\par
Assume that $C=C_1$.
Then we have $n=1,3$ by condition
(4). If
$n=1$ then there is a
 2-complement of $K_S+C+B$ as before.
Consider the case $n=3$. Then
$B^{\mt{hor}}=\frac12B_1+\frac12B_2$ or
$B^{\mt{hor}}=\frac12B_1$. The divisor $K_S+C+B$ doesn't have 1-,2-,3-,4-
and 6-complement if and only if the divisor $K_C+\Diff_C(B)$ doesn't have
1-,2-,4- and 6-complement, that is, (after simple calculations)
\begin{equation}\label{eq1}
\tag{*}
(C,\Diff_C(B))=(\PP^1,(3/5+\varepsilon_1)P_1+(2/3+\varepsilon_2)P_2+
(5/7+\varepsilon_3)P_3),
\end{equation}
where $\varepsilon_1+\varepsilon_2+\varepsilon_3<\frac2{105}$ and
$\varepsilon_i\ge 0$ for all $i$. By lemma \ref{compl4} the divisor
$K_S+C+B$ is 12-complementary (the index 12 is not always a minimal one).
\par
Let $\nu\colon S\stackrel{\psi_1}{\longrightarrow}
S'\longrightarrow Z$ be a contraction of all curves in the fibres of $\nu$ (with the
help of log minimal model program) with $(K_S+C+\down{13B}/12)\cdot E>0$. Since
$(K_S+C+\down{13B}/12)\cdot C=0$ then $\psi_1$
doesn't contract the curves intersecting $C$.
We get that the divisor
$(K_{S'}+C'+\down{13B'}/12)$ is nef and in particular, it is nef over
$Z$, where $C'$ and $B'$ are the images of $C$ and $B$.
The cone $\NE(S/Z)$ is polyhedral and generated by contractible extremal curves
\cite[proposition 2.5]{Sh2}. Let
$\nu\colon S\stackrel{\psi_1}{\longrightarrow} S'
\stackrel{\psi_2}{\longrightarrow} \overline S\longrightarrow Z$
be a contraction of all curves not intersecting $C'$ in the fibres of
$\nu$.
Then $\rho(\overline S/Z)=1$. The pair $(\overline C,\Diff_{\overline
C}(\overline B))$ is the same one as in (\ref{eq1}), where
$\overline C$ and $\overline B$ are the images of $C'$ and $B'$. Note that either
$\Diff_{\overline C}(0)=0$, or $\Diff_{\overline
C}(0)=\frac23P_2$ and $\varepsilon_2=0$.
\par
1). Consider the case $\Diff_{\overline C}(0)=0$. Since
$\overline C\cdot\overline B^{\mt{hor}}=0$ then $\overline{C}^2\le 0$.
Hence, the linear system $|\overline C+mf|$ gives a birational
morphism $\psi\colon\overline S\to \FFF_k$, where
 $f$ is a general fiber and $m\gg 0$
\cite[proposition 1.10]{KudLd}. We obtain that
$\overline S\cong \FFF_k$.
\par
2). Consider the case $\Diff_{\overline C}(0)=\frac23P_2$ and
$\varepsilon_2=0$. Let $\phi\colon \overline S'\to \overline S$
be the blow-up with the unique exceptional curve at the point
$P_2$ such that
$\Sing \overline S'\cap \overline C'=\emptyset$ and the divisor
$K_{\overline S'}$ is $\phi$-nef. Put $f_2=\Exc
\phi$. By the same argument as in the previous case  the linear system $|\overline
C'+mf|$ gives a birational morphism $\psi\colon\overline S'\to
\FFF_k$. Let $E=\Exc \psi$, $\widehat D=\psi(\overline
C')+\down{13\psi(\overline B')}/12$ and $\psi(E)=P$. Then either
$$
\big(\FFF_k\ni P,\widehat D\big)\an \big(\CC^2_{x,y}\ni 0,
(2/3)\{y=0\}+ (1/2)\{x(x+y^t)=0\}\big),
$$
or
$$
\big(\FFF_k\ni P,\widehat D\big)\an \big(\CC^2_{x,y}\ni 0,
(2/3)\{y=0\}+ (1/2)\{x^2+y^t)=0\}\big).
$$
Since $a(E,\widehat D)=0$ and $f_2^2<-1$ then $\psi$ is a weighted blow-up
with weights (1,3) or (2,3) (cf. \cite[lemma 5.5]{KudLd}).
In the second case if we take a blow-up with weights
(2,3) then the following condition must be satisfied: $t\ge 2$.
The result is summarized
in the next example.
\begin{example}\label{mainind1}
1). Let
$$\big(\FFF_k,D^+\big)=\big(\FFF_k,E_{\infty}+(1/2)E_1+(1/2)E_2+(7/12)f_1+(2/3)f_2+
(3/4)f_3\big),$$ where $E_{\infty}$ is a minimal section, $E_i$ is a zero section,
$f_i$ is a fiber. Let $h\colon S\to \FFF_k$ be a birational contraction:
\begin{gather*}
K_S+\widetilde E_{\infty}+(1/2)\widetilde E_1+(1/2)\widetilde E_2+
(7/12)\widetilde f_1+(2/3)\widetilde f_2+(3/4)\widetilde f_3+\\+\sum a_i\Delta_i
=h^*(K_{\FFF_k}+D^+)\equiv 0.
\end{gather*}
Assume that $a_i\in\{0\}\cup[1/2,1)\cup\ZZ/(12)$,
$h(\Delta_i)\not\in E_{\infty}$ for all $i$. Since the pair
$(\FFF_k,D^+)$ is kawamata log terminal outside $E_{\infty}$ then
there is only finite number of such surfaces $S$ by lemma
3.1.9 \cite{PrLect}.
The log surface
\begin{gather*}
\big(S,D\big)=\big(S,\widetilde E_{\infty}+(1/2)\widetilde
E_1+(1/2)\widetilde E_2+ (3/5+\varepsilon_1)\widetilde
f_1+(2/3+\varepsilon_2)\widetilde f_2+\\+
(5/7+\varepsilon_3)\widetilde f_3+\sum_{i:\ a_i\ge 1/2}
a_i\Delta_i \big)
\end{gather*}
satisfies the condition of theorem \ref{mainind}, where
$\varepsilon_1+\varepsilon_2+\varepsilon_3\le\frac2{105}$ and
$\varepsilon_i\ge 0$ for all $i$. If $5/7+\varepsilon_3<8/11$
then we have 10-complement of $K_S+D$. If $5/7+\varepsilon_3\ge 8/11$
then we have 12-complement. Note also that sometimes we can easily
change the coefficients
$a_i$, and we can contract $\widetilde E_{\infty}$ if $k>0$.
\par
2). Let $\overline S\to \PP^1$ be an extremal generically
$\PP^1$-fibration, that is, $\rho(\overline S/\PP^1)$=1. Assume that
$\Sing \overline{S}\subset f_2$, where the fiber $f_2$ is shown in one of
the following figures.
\\
\begin{center}
\begin{picture}(200,70)(0,-20)
\put(10,20){\fbox{$\overline E_{\infty}$}}
\put(35,24){\line(1,0){20}} \put(59,24){\circle*{8}}
\put(55,30){\scriptsize{-2}} \put(63,24){\line(1,0){20}}
\put(87,24){\circle*{8}} \put(83,30){\scriptsize{-2}}
\put(91,24){\line(1,0){20}} \put(119,24){\circle{16}}
\put(115,21){\footnotesize{$f_2$}} \put(127,24){\line(1,0){20}}
\put(151,24){\circle*{8}} \put(147,30){\scriptsize{-3}}
\put(186,40){\fbox{$\overline B_1$}} \put(151,24){\line(2,1){35}}
\put(151,24){\line(2,-1){35}} \put(186,0){\fbox{$\overline B_2$}}
\put(95,-15){Fig. 1}
\end{picture}
\end{center}

\begin{center}
\begin{picture}(200,90)(0,-20)
\put(10,20){\fbox{$\overline E_{\infty}$}}
\put(35,24){\line(1,0){20}} \put(59,24){\circle*{8}}
\put(55,30){\scriptsize{-3}} \put(63,24){\line(1,0){20}}
\put(91,24){\circle{16}} \put(87,21){\footnotesize{$f_2$}}
\put(99,24){\line(1,0){20}} \put(123,24){\circle*{8}}
\put(119,30){\scriptsize{-2}} \put(127,24){\line(1,0){20}}
\put(151,24){\circle*{8}} \put(147,30){\scriptsize{-2}}
\put(186,40){\fbox{$\overline B_1$}} \put(151,24){\line(2,1){35}}
\put(151,24){\line(2,-1){35}} \put(186,0){\fbox{$\overline B_2$}}
\put(95,-15){Fig. 2}
\end{picture}
\end{center}

\begin{center}
\begin{picture}(200,90)(0,-40)
\put(10,20){\fbox{$\overline E_{\infty}$}}
\put(35,24){\line(1,0){20}} \put(59,24){\circle*{8}}
\put(55,30){\scriptsize{-3}} \put(63,24){\line(1,0){20}}
\put(91,24){\circle{16}} \put(87,21){\footnotesize{$f_2$}}
\put(99,24){\line(1,0){20}} \put(123,24){\circle*{8}}
\put(119,30){\scriptsize{-2}} \put(127,24){\line(1,0){20}}
\put(151,24){\circle*{8}} \put(147,30){\scriptsize{-2}}
\put(114,-10){\fbox{$\overline B_1$}}
\put(123,24){\line(0,-1){20}} \put(95,-40){Fig. 3}
\end{picture}
\end{center}

Let us consider the minimal resolution of a surface
$\overline S$ and contract all
(--1) curves not intersecting the proper transform of
$\overline E_{\infty}$. We obtain $\FFF_k$. On the ruled surface
$\FFF_k$ the image of $\overline E_{\infty}$ is a minimal section,
the images of $\overline B_1$ and $\overline B_2$ from figures
1 and 2 are the sections, the image of $\overline B_1$ from figure
3 is a 2-multi-section.
\par
Consider the log surface
$$\big(\overline S,D^+\big)=\big(\overline S,\overline E_{\infty}+\overline B+
(7/12)f_1+(3/4)f_3\big),$$ where $\overline B=\frac12\overline
B_1+\frac12\overline B_2$ (in the case of figure 1 or 2) or
$\overline B=\frac12\overline B_1$ (in the case of figure 3).
Arguing as above in the previous point of example we
can construct the birational morphisms
$h\colon S\to \overline S$.
It is clear that the same statements take place about the structure of
$h$ and the complements of $K_S+D$.
\end{example}

\par C). Let $C$ be in a fiber of $\nu$. Put $P=\nu(C)$ and $f=\nu^{-1}(P)$.
The case, where the general fiber is an elliptic curve is considered in
III, B). Therefore we assume that the general fiber of
$\nu$ is a rational curve.
The divisor $K_C+\Diff_C(B)$ has an
$n$-semicomplement of minimal index, where $n\in\{1,2,3,4,6\}$
(see proposition \ref{compl1}).
Assume that there is a horizontal component
$B_i$ with a coefficient
$b_i\not\in\ZZ/(n+1)$.
Then considering the
divisor
$K_S+C+B-\varepsilon B_i$ we reduce our problem to case I.
Therefore we assume that all horizontal components of
$B$ have the coefficients from the set
$\ZZ/(n+1)$. Let us show that this possibility is
impossible.
\par
Let $\nu\colon S\stackrel{\psi}{\longrightarrow}
S'\longrightarrow Z$ be a contraction of all curves in the fiber $f$ not lying in
$C$. Put $C'=\psi(C)$, $B'=\psi(B)$. Since
$K_S+C+B\equiv 0$ over $Z$ then
$(C,\Diff_C(B))=(C',\Diff_{C'}(B'))$.
\par
Assume that  $C$ is a reducible curve. Then $n=1,2$. If
$n=1$ then all coefficients of horizontal components of $B$
are equal to 1/2 by condition (4).
Since the divisor $K_{S'}+C'+B'$ is divisorial log terminal
and numerical trivial over $Z$ then there is a divisor in
$\Diff_{C'}(B')$ with a coefficient 1/2.
A contradiction with $n=1$. If $n=2$
then we have same contradiction.
\par
The case, where $C$ is an irreducible curve, is considered similarly.
\par {\bf Case III.} Assume that $K_S+C+B\equiv 0$.
\begin{definition}
Let $D$ be a $\QQ$-divisor on a projective variety
$X$. Define {\it the numerical dimension} of a divisor $D$:
$$
\num(D)=\max\{\num(D')\mid D'\ge 0,\ \Supp D'\subseteq \Supp D\}.
$$
{\it The linear dimension $($Iitaka dimension$)$} $\kappa(D)$
is defined similarly.
\end{definition}

We have $\kappa(B)=\num(B)$ \cite[proposition 2.12]{Sh2}.
Let us consider all possibilities of $\num(B)$
case by case.
\par A). Let $\num(B)=2$. The cone $\NE(S)$ is polyhedral and generated by
contractible extremal curves since there exists a divisor
$B'$ ($\Supp B'\subseteq\Supp B$) such that the divisor
$-(K_S+C+B)+\varepsilon B'$ is nef and big
\cite[proposition 2.5]{Sh2}. Let the divisor
$K_C+\Diff_C(B)$ be $n$-semicomplementary. Let us contract all exceptional curves
$E$ with
$(K_S+C+\down{(n+1)B}/n)\cdot E>0$. We obtain either $\nu\colon
S\stackrel{\psi}{\longrightarrow} S'
\stackrel{\nu'}{\longrightarrow} Z\cong \PP^1$ and
$(K_{S'}+C'+\down{(n+1)B'}/n)\cdot f>0$, where $f$ is a general fiber of
$\nu'$, or $\psi\colon S\to S'$ and
$-(K_{S'}+C'+\down{(n+1)B'}/n)$ is a nef divisor.
By the construction none component of
$C$ is contracted by $\psi$, and $C'$ doesn't lie in the fibres of
$\nu'$ in the first case. Therefore
$C$ has the horizontal components of
$\nu$ in the first case. Moreover,
since $\num(B)=2$ then we have $B^{\mt{hor}}=\frac12B_1+\frac12B_2$ or
$B^{\mt{hor}}=\frac12B_1$ by condition (4).
Thus the first case is reduced to case
II, B). This new possibility was included in example
\ref{mainind1}.
\par
In the second case the divisor $K_{S'}+C'+\down{(n+1)B'}/n$ is
$n$-complementary without
condition (4) on the coefficients of a boundary $B$
\cite{Sh2}
(here it is essential that the cone
$\NE(S')$ is polyhedral and generated by contractible extremal curves).
By proposition \ref{compl3} the divisor $K_S+C+B$ is $n$-complementary.
\par B). Let $\num(B)=1$.
Then for some divisor $B'$ ($\Supp B'\subseteq \Supp B$) the linear system
$|B'|$ gives a fibration $\nu\colon S\to Y$ and the divisor $B$ lies in the fibres.
If some component of $C$ is a section then the general fiber is
$\PP^1$. Hence there is a multi-section
$C_1\subset C$. By lemma \ref{compl4} our theorem is proved. Therefore we
may assume that $C$ lies in the fibres and the general fiber is an elliptic curve.
Arguing as above, we contract all curves
$E\not\subset\Supp C$ such that $(K_S+C+\down{(n+1)B}/n)\cdot
E\ge 0$, where $n$ is the semicomplementary index of
$K_S+\Diff_C(B)$.
As a result we obtain
$\nu\colon S\stackrel{\phi}{\longrightarrow} S'
\stackrel{\nu'}{\longrightarrow} Y$, where all fibers of $\nu'$ are irreducible,
except the fiber consisting of $C'$. Let
$\psi\colon S'\to S''$ be the contraction of components of $C'$. We get the model with
$\rho(S'')=2$. The cone of surface $S''$ has two extremal rays:
a fiber of $\nu'$, a (multi-)section $E$. If
$(K_{S''}+C''+\down{(n+1)B''}/n)\cdot E\le 0$ then the divisor
$-(K_{S'}+C'+\down{(n+1)B'}/n)$ is nef. In this case the theorem is proved
in \cite{Sh2} without
condition (4) on the coefficients of a boundary
$B$.
If $(K_{S''}+C''+\down{(n+1)B''}/n)\cdot E>0$ then we have a fibration
$\theta\colon S\to Z\cong \PP^1$ and $C$ contains a section of
$\theta$. The horizontal part of divisor
$B$ for $\theta$ is $B^{\mt{hor}}=\frac12B_1+\frac12B_2$ or
$B^{\mt{hor}}=\frac12B_1$ by condition
(4). Since $K_S+C+B\equiv 0$ over $Y$ and
$K_{S''}+C''\equiv 0$ over $Y$ then $\Diff_C(B)=\Diff_{C''}(0)$.
Hence the coefficients of $\Diff_C(B)$ are standard. In particular,
if $n=1,3$ then we have $\Diff_C(B)=P_1+P_2$ and
$\Diff_C(B)=\frac23P_1+\frac23P_2+\frac23P_3$ respectively.
Therefore, if $n=1,2,3,4,6$ then the divisor $K_S+C+B$ is
2-,2-,6-,4-,6-complementary respectively by lemma \ref{compl4}
\par C). Let $\num(B)=0$, that is,  the divisor $B$ is contracted.
Then our theorem is true without condition
(4) on the coefficients of a boundary
$D$ \cite{Sh2}.
\end{proof}
\end{theorem}

\begin{remark}
Condition (3) of theorem \ref{mainind} can be replaced by
one of the following more strong conditions
(see \cite[proposition 2.5]{Sh2}):
\begin{enumerate}
\item[$(3^1)$] the divisor $-(K_S+D)$ is big;
\item[$(3^2)$] the cone
$\NE(S)$ is polyhedral and generated by contractible extremal curves;
\item[$(3^3)$] the divisor $-(K_S+D)$ is semi-ample;
\item[$(3^4)$] there exists a contraction $\nu\colon S\to Z$ with the following
property:\\ if $(K+D)\cdot E=0$ then $E\subset \Exc\nu$.
\end{enumerate}
In example \ref{mainind1} the log surfaces $(S,D)$ satisfy the
conditions enumerated, of course except condition $(3^1)$.
\end{remark}

The next corollary is very important for the applications.
\begin{corollary}\label{cor1}
Let $(S,D=\sum d_iD_i)$ be a projective log surface.
Assume that
\begin{enumerate}
\item the divisor $K_S+D$ is kawamata log terminal;
\item the divisor $-(K_S+D)$ is nef;
\item there exists a $\QQ$-complement of $K_S+D$;
\item $d_i\ge \frac12$ for all $i$;
\item there exists an effective
$\QQ$-divisor $D'\ge D$ such that the divisor $-(K_S+D')$ is nef
and the pair $(S,D')$ is not kawamata log terminal.
\end{enumerate}

Then there is $1$-, $2$-, $3$-, $4$- or $6$-complement of $K_S+D$ which is not
kawamata log terminal, except the cases appearing in example
$\ref{mainind1}$.
\par
Besides, if there is an infinite number of divisors
$E$ with a discrepancy
$a(E,D)=-1$ and the pair $(S,D')$ is log canonical
then there is $1$- or $2$-complement of
$K_S+D$ which is not kawamata log terminal.

\begin{proof} Replacing the divisor $D'$ with suitable $D+\lambda(D'-D)$, where
$\lambda>0$ we may assume that the divisor $K_S+D'$ is log canonical but not
kawamata log terminal.
\par
At first let us prove that there is a $\QQ$-complement of
$K_S+D'$. If the divisor
$-(K_S+D)$ is big then the cone
$\NE(S)$ is polyhedral and generated by contractible extremal curves
and we obtain the required statement \cite[proposition 2.5]{Sh2}.
Therefore it can be assumed that
the linear system
$|-m(K_S+D)|$ gives a fibration $\nu\colon S\to Z$, where $m\gg 0$.
Adding the required number of general fibres of $\nu$ to
the divisor $K_S+D'$ we have our statement.
\par
Let $f\colon \widetilde S\to S$ be a minimal log terminal modification of the pair
$(S,D')$ \cite[definition 3.1.3]{PrLect}. We have
$$
K_{\widetilde S}+\sum E_i+\down{\widetilde D'}+\fr{\widetilde D'}=
f^*(K_S+D'),
$$
where $\widetilde D'$ is a proper transform of $D'$. Put
$$
\widetilde D=\sum E_i+\down{\widetilde D'}+\sum_{i\colon \Supp
\widetilde D_i \not\subset \Supp\down{\widetilde D'}}
d_i\widetilde D_i,
$$
where $\widetilde D_i$ is a proper transform of $D_i$.
Thus, the statement of corollary must be proved for the divisor
$K_{\widetilde S}+\widetilde D$. If the divisor $-(K_{\widetilde
S}+\widetilde D)$ is nef then it is nothing to be proved by theorem
\ref{mainind}. Therefore it can be assumed that the divisor
$-(K_{\widetilde S}+\widetilde D)$ is not nef.
A $\QQ$-complement of a divisor $K_{\widetilde S}+\widetilde D$
is denoted by $\widetilde\Theta$. We can assume that
$\down{\widetilde\Theta}=\down{\widetilde D}$.
\par
Let us prove that we can contract all exceptional curves
$E$ such that $(K_{\widetilde S}+\widetilde D)\cdot E>0$ on every step.
\par A). Assume that $\nu_{\mt{num}}(\fr{\widetilde\Theta})=2$.
Then arguing as in the proof of theorem
\ref{mainind} (case III,A)) the cone
$\NE(\widetilde S)$ is polyhedral and generated by contractible extremal curves.
Q.E.D.

\par B). Assume that $\nu_{\mt{num}}(\fr{\widetilde\Theta})=1$. By
proposition 2.12 \cite{Sh2} for some divisor
$\widetilde\Theta'$ ($\Supp\widetilde\Theta'\subset
\Supp\fr{\widetilde\Theta$}) the linear system
$|\widetilde\Theta'|$ gives a fibration $\nu\colon \widetilde S\to
Z$ and a divisor $\fr{\widetilde\Theta}$ lies in the fibres of  $\nu$.
If $(K_{\widetilde S}+\widetilde D)\cdot E>0$ then a curve $E$ lies
in the fibres of $\nu$. Therefore it can be contracted.
\par C). Assume that $\nu_{\mt{num}}(\fr{\widetilde\Theta})=0$. Then
a divisor $\fr{\widetilde\Theta}$ is contracted by the definition.
\par
Thus we get a birational morphism $\phi\colon \widetilde S\to
\overline S$. It is clear that $\phi$ doesn't contract the components of
$\down{\widetilde D}$, and the curve contracted intersects some component
$\widetilde\Theta_1$ of
$\down{\widetilde D}$ on every step. Put $\overline D=\phi(\widetilde D)$.
\par
It remains to prove that
an $n$-complement $\overline D^+$ of
$K_{\overline S}+\overline D$ induces an $n$-complement of
$K_{\widetilde S}+\widetilde D$ ($n=$1,2,3,4 or 6). Put
$$
K_{\widetilde S}+\widetilde D^+=\phi^*(K_{\overline S}+\overline
D^+).
$$
We must prove that
\begin{equation*}\tag{**}\label{**}
\down{\widetilde D}+\frac{\down{(n+1)\fr{\widetilde D}}}n\le
\down{\widetilde D^+}+\fr{\widetilde D^+}.
\end{equation*}

By the above this requirement is enough to check in the case, when
$\phi$ is a contraction of the unique curve $E$.
Let $P=\phi(E)$.
By the classification of
two-dimensional log terminal pairs
\cite[theorem 2.1.2]{PrLect} and by condition (4)
we conclude that there are at most one divisor of
$\fr{\overline D}$ passing through the point $P$
and $(\overline S\ni P)$ is a cyclic singularity.
Let $d_1\overline D_1$ be a divisor passing through the point $P$.
If the coefficient of divisor
$\overline D_1$ in $\overline D^+$ is more then
$d_1$ then we consider it instead of
$d_1$.
Since the divisor
$K_{\down{\overline D}}+\Diff_{\down{\overline D}} (\fr{\overline D})$ is
$n$-semicomplementary there are the following cases (the case $n=1$ is obvious).
\par 1). $(\overline S\ni P, \overline\Theta_1+d_1\overline D_1)\an
(\CC^2\ni 0, \{x=0\}+d_1\{y=0\})/\ZZ_2(1,1)$ and $n=4$. Then by
proposition \ref{compl3} requirement (\ref{**})
must be checked for
$d_1\in(\frac12,\frac35)$.
\par 2). $(\overline S\ni P, \overline\Theta_1+d_1\overline D_1)\an
(\CC^2\ni 0, \{x=0\}+d_1\{y=0\})/\ZZ_2(1,1)$ and $n=6$. Then by
proposition \ref{compl3} requirement (\ref{**}) must be checked for
$d_1\in(\frac23,\frac57)$.
\par 3). $(\overline S\ni P, \overline\Theta_1+d_1\overline D_1)\an
(\CC^2\ni 0, \{x=0\}+d_1\{y=0\})/\ZZ_3(q,1)$ and $n=6$. Then by proposition
\ref{compl3} requirement (\ref{**}) must be checked for
$d_1\in(\frac12,\frac47)$.
\par 4). $(\overline S\ni P, \overline\Theta_1+d_1\overline D_1)\an
(\CC^2\ni 0, \{x=0\}+d_1\{y=0\})$.\\
Requirement (\ref{**}) is equivalent to the following one:
\begin{equation*}\tag{***}
-\frac{\down{(n+1)a(E,\overline\Theta_1+d_1\overline D_1)}}n\le
-a\Big(E,\overline\Theta_1+\frac{\down{(n+1)d_1}}n\overline
D_1\Big).
\end{equation*}
Since $a(E,\overline\Theta_1+d_1\overline D_1)\le
-1/2$ then $\phi$ is a toric blow-up.
Requirement (***) in cases 1),2) and 3) is checked directly.
In case 4) the weights of weighted blow-up $\phi$ are denoted by
$(\alpha,\beta)$. Then either $(\alpha,\beta)=(\alpha,1)$ and $d_1\ge
1/2$, or $(\alpha,\beta)=(\alpha,2)$, $d_1\ge 3/4$ and $n=4$, or
$(\alpha,\beta)=(\alpha,3)$, $d_1\ge 5/6$ and $n=6$.
Now requirement (***) is also fulfilled by
direct calculation.
\end{proof}
\end{corollary}

\begin{corollary} Under the notation of corollary
$\ref{cor1}$ let us decline condition
$(4)$ on a boundary
$D$. Then there is $n$-complement of
$K_S+D$ which is not kawamata log terminal, where
$n=1$, $2$, $3$, $4$, $5$, $6$, $7$,
$8$, $9$, $10$, $11$, $12$, $13$, $14$, $15$, $16$, $17$, $18$,
$19$, $20$, $21$, $22$, $23$, $24$, $25$, $26$, $27$, $28$, $29$,
$30$, $31$, $35$, $36$, $40$, $41$, $42$, $43$, $56$ or $57$.
\par
Besides, if there is an infinite number of divisors
$E$ with a discrepancy
$a(E,D)=-1$ and the pair $(S,D')$ is log canonical
then there is $1$-, $2$- or $6$-complement of
$K_S+D$ which is not kawamata log terminal.
\begin{proof}
By the proof of corollary \ref{cor1} it follows that the divisor
$K_S+D'$ has a $\QQ$-complement. Therefore
we obtain our statement by theorem 2.3 \cite{Sh2}.
\end{proof}
\end{corollary}

\section{\bf {Exceptional non-rational log surfaces}}

\begin{theorem}\label{main2}
$($cf. \cite[proposition 9.2.2]{PrLect}$)$ Let $(S,D=\sum
d_iD_i)$ be a projective log surface, where $D$ is a boundary.
Assume that the following conditions are satisfied:
\begin{enumerate}
\item there exists a $\QQ$-complement $\Theta=\sum\theta_i \Theta_i $
of $K_S+D$;
\item the surface $S$ is non-rational and the pair $(S,D)$
is exceptional;
\item we have $\Theta\ne 0$ or $S$ has a non Du Val singularity.
\end{enumerate}

Then there is $2$-, $3$-, $4$- or $6$-complement of
$K_S+D$.
Besides, one of the following cases takes place.
\begin{enumerate}
\item[1).] $S\cong C\times\PP^1$, where $C$ is an elliptic curve,
$D_i$ are the sections of corresponding $\PP^1$-bundle.
\item[2).] $S\cong \PP_C(\mathcal E)$, where $\mathcal E$ is an indecomposable
vector bundle of degree
$1$ on an elliptic curve $C$.
Up to multiplication by an invertible sheaf, $\mathcal E$ is a nontrivial
extension
$$
0\longrightarrow\OO_C\longrightarrow\mathcal
E\longrightarrow\OO_C(O)\longrightarrow 0.
$$
Then $D_i\sim 2E-f^*(O+t_i)$ or $D_i\sim 4E-2f^*O$, where $E$ is a section of
$f\colon S\to C$ and $t_i$ is an element of order $2$ in
$\Pic(C)$.
\end{enumerate}
\begin{proof}
Let $\phi\colon\widetilde S\to S$ be a minimal resolution.
Then $K_{\widetilde S}+\widetilde\Theta=\phi^*(K_S+\Theta)\equiv
0$. Condition (3) implies $\widetilde\Theta\ne 0$, that is,
$\kappa(\widetilde S)=-\infty$. Let $S_{\min}$ be a minimal model of
$\widetilde S$. By the condition $S_{\min}$ is a minimal ruled surface over
a curve $C$ with $p_a(C)\ge 1$.
The image of divisor $\widetilde\Theta$ on $S_{\min}$ is denoted by
$\overline\Theta$. If there is an irreducible curve
$E$ with $E^2<0$ on $S_{\min}$ then
$$
(K_{S_{\min}}+E)\cdot E=\Big(-\sum_{i\colon\overline\Theta_i\ne E}
\theta_i\overline\Theta_i\Big)\cdot E+\varepsilon E^2<0,
$$
where $\varepsilon>0$. Hence $p_a(E)=0$ and we have a contradiction with
$p_a(C)\ge 1$.
\par
Since $0\ge 8-8p_a(C)=K^2_{S_{\min}}=\overline\Theta^2\ge 0$
then $p_a(C)=1$,
$\overline\Theta_i^2=\overline\Theta_i\cdot\overline\Theta_j=0$
for all $i,j$. Since $\theta_i<1$ for all $i$ then the pair
$(S_{\min},\overline\Theta)$ is terminal. Therefore
$\PP_C(\mathcal E)\cong S_{\min}\cong\widetilde S\cong S$, where
$\deg \mathcal E\ge 0$. By chapter 5 \cite{Hart} and by examples 1.1,
2.1 \cite{Sh2} we obtain the remaining statements.
\end{proof}
\end{theorem}

\begin{remark} \cite[example 2.1]{Sh2} In case 2) of theorem \ref{main2}
the linear system $|4E-2f^*O|$ gives a structure of
elliptic fibration with three degenerate (double) fibres, which are
linear equivalent to
$2E-f^*(O+t_i)$.
\end{remark}

\begin{corollary} Under the conditions of theorem $\ref{main2}$
we have $\delta(S,D)=0$.
\end{corollary}

\section{\bf {Construction of models of log del Pezzo surfaces with $\delta=0$}}

The classification of exceptional log surfaces with
$\delta=1,2$ was given in the papers
\cite{KudLd}, \cite{Sh2}. The exceptional non-rational log surfaces were completely
classified in theorem \ref{main2}.
Thus it remains to study the last remaining case
-- the exceptional rational log surfaces with $\delta=0$.
\begin{definition}
The pair $(S,D)$ is called {\it a log del Pezzo surface}, where $D$ is a boundary,
if the following conditions are satisfied:
\begin{enumerate}
\item the divisor $-(K_S+D)$ is nef;
\item the divisor $K_S+D$ is log canonical;
\item there exists a $\QQ$-complement of
$K_S+D$.
\end{enumerate}
\end{definition}

Let us consider the limiting case.
\begin{example}\label{ex1}\cite[remark 1.2]{Bl}, \cite[examples 4.2, 5.3]{Zh1}
1). Let $\overline S=C\times C$, where
$C=\CC/(\ZZ+\varepsilon_3\ZZ)$ is an elliptic curve and
$\varepsilon_3$ is a primitive root of  unity of order 3. The group
$\ZZ_3$ acts on the curve $C$ by the multiplication on $\varepsilon_3$.
Then $S=\overline S/\ZZ_3$ is a surface with $3K_S\sim 0$,
$\rho(S)=4$, and $\Sing S$ consists of nine singularities  $\frac13(1,1)$.
\par
2). Let the surface $\overline S=J(C)$ be the jacobian of hyperelliptic curve
$C\colon y^2=x^5-1$ of genus 2. The group
$\ZZ_5$ is generated by the automorphism $(x,y)\longmapsto
(\varepsilon_5x,y)$ of curve $C$, where $\varepsilon_5$ is a
primitive root of  unity of order 5.
Then $S=\overline S/\ZZ_5$ is a surface with
$5K_S\sim 0$, $\rho(S)=2$, and $\Sing S$
consists of five singularities $\frac15(2,1)$.
\par
3). Let us consider three irreducible
curves $\overline
E_1\sim\OO_{\PP^2}(1)$, $\overline E_2\sim\overline
E_3\sim\OO_{\PP^2}(4)$ on
$\PP^2$.
The curve $\overline E_2$ and the curve $\overline E_3$
has three ordinary double points
$Q_5$, $Q_6$, $Q_7$ and
$Q_8$, $Q_9$, $Q_{10}$ respectively. The line $\overline E_1$
intersects the curves $\overline E_2$ and $\overline E_3$ at the points
$Q_1$, $Q_2$, $Q_3$ and $Q_4$. The curve $\overline E_2$ intersects
$\overline E_3$ at the points $Q_1$,\ldots,$Q_{10}$. Let us take
the usual blow-ups of $\PP^2$
at the points
$Q_1$,\ldots,$Q_{10}$ and contract the proper transforms of the curves
$\overline E_1$, $\overline E_2$, $\overline E_3$. We get a surface
$\widetilde S$ with $3K_{\widetilde S}\sim 0$,
$\rho(\widetilde S)=8$, and $\Sing \widetilde S$ consists of three singularities
$\frac13(1,1)$. When contracting $(-2)$ curves on $\widetilde
S$ we get a surface $S$ with $3K_S\sim 0$ and with three non Du Val singularities
$\frac13(1,1)$.
\end{example}

\begin{theorem}\label{main3}
Let $S$ be a rational exceptional log del Pezzo surface
$(D=0)$. Assume that $a(E,0)>-1/2$ for all $E$ and there is no
a $\QQ$-complement $\Theta=\sum\frac12\Theta_i$ of
$K_S$. Then the surface $S$ is of example $\ref{ex1}$.
\begin{proof} The surface  $S$ must have a non Du Val singularity otherwise
$h^0(S,\OO_S(-K_S))\ge 1$.
\begin{lemma}\label{lem1}
Let $(X\ni P)$ be a two-dimensional non Du Val singularity.
Assume that
\begin{eqnarray*}
M=\min\{ a(E,0)\mid E\  \text{is an exceptional divisor}\} >-\frac12.
\end{eqnarray*}
Then $(X\ni P)\an (\CC^2\ni 0)/\ZZ_{2n+1}(n,1)$, where $n\ge 1$. In particular,
$M=-\frac{n}{2n+1}\le-\frac13$.
\begin{proof}
If $(X\ni P)$ is a non-cyclic singularity then the blow-up of central vertex of
minimal resolution graph gives a discrepancy
$\le-\frac12$. Therefore $(X\ni P)$ is a cyclic singularity.
Let $\widetilde X\to X$ be a blow-up with the unique exceptional exceptional curve
$E$ such that its self-intersection index
$k$ on the minimal resolution
of $(X\ni P)$ is at most $-3$. Then
$$
-\frac12<a(E,0)=-1+\frac{-2+\Diff_E(0)}{E^2}=-1+\frac{-2+\frac{m_1-1}{m_1}+
\frac{m_2-1}{m_2}}{-k+\frac{q_1}{m_1}+ \frac{q_2}{m_2}}.
$$
Hence $(m_1,q_1)=(1,0)$, $(m_2,q_2)=(m_2,m_2-1)$, $k=3$.
\end{proof}
\end{lemma}

Let $P_1,\ldots,P_r$ be non Du Val singularities of $S$ of types
$\frac1{2n_1+1}(n_1,1)$,\ldots, $\frac1{2n_r+1}(n_r,1)$ respectively. Let
$f\colon \widetilde S\to S$ be a minimal resolution. Then
$K_{\widetilde S}+\Delta=f^*K_S$. By lemma \ref{lem1}
$$
h^2(\widetilde S,\OO_{\widetilde S}(-2K_{\widetilde
S}-\down{3\Delta}))= h^0(\widetilde S,\OO_{\widetilde
S}(3K_{\widetilde S}+\down{3\Delta}))=0,
$$
except the case $n_1=\ldots=n_r=1$ and $K_S\equiv 0$.
Let us determine the remaining possibilities of
$n_1,\ldots,n_r$. By Riemann-Roch theorem and Noether's formula
we have the next system
$$
\left\{
\begin{array}{l}
0=h^0(\widetilde S,\OO_{\widetilde S}(-2K_{\widetilde
S}-\down{3\Delta}))=3K_S^2
+r+1-3\cdot\sum_{i=1}^r\frac{n_i}{2n_i+1}+\\
\hspace{6.5cm}
+h^1(\widetilde S,\OO_{\widetilde S}(-2K_{\widetilde S}-\down{3\Delta}))\\
K_S^2-\sum_{i=1}^r\frac{n_i}{2n_i+1}+\rho(S)+n_1+\ldots+n_r\le 10.  \\
\end{array}
\right.
$$

\begin{lemma}\label{lem2}\cite[corollary 9.2]{KeM}
Let $X$ be a rational surface with kawamata log terminal singularities
and with $\rho(X)=1$. Then
$$\sum_{P\in\Sing X}\frac{m_P-1}{m_P}\le 3,$$ where
$m_P$ is the order of the local fundamental group
$\pi_1(U_P\setminus\{P\})$ $(U_P$ is a sufficiently small neighborhood
of $P)$.
\end{lemma}

Taking into account $K_S^2\ge 0$ and lemma \ref{lem2} we obtain
$K_S^2=0$, $\rho(S)=2$, $n_1=\ldots=n_5=2$, $\Sing
S=\{P_1,\ldots,P_5\}$ by the system. Moreover, we have $K_S\equiv 0$.
Indeed, let $K_S\not\equiv 0$. If there is a curve $E$ with $E^2<0$ on
$S$ then we can contract it and obtain a contradiction with lemma
\ref{lem2}. Therefore we have generically $\PP^1$-fibration $S\to Z$, but it is impossible
by classification of such fibrations \cite[theorem 7.1.12]{PrLect}.
\par
Let $\overline S\to S$ be a canonical cover. There are two cases
\cite[theorem C]{Bl}.
\par A). Let $\overline S$ be an abelian surface.
Then theorem C \cite{Bl} implies that $S$ is a surface from example
\ref{ex1} (points 1) or 2)) and
$n_1=\ldots=n_9=1$ or $n_1=\ldots=n_5=2$ respectively.
\par B). Let $\overline S$ be a $K3$-surface.
Then by theorem 5.1 \cite{Zh1} we get that $n_1=n_2=n_3=1$.
It is clear that one of the minimal models of
$\widetilde S$ is
$S_{\min}\cong \PP^2$. Let $\varphi\colon\widetilde S\to
S_{\min}$ be a corresponding birational morphism. Put
$\overline\Delta=\varphi(\Delta)=\frac13(\varphi(E_1)+\varphi(E_2)+\varphi(E_3))$,
where $E_i$ is an exceptional curve over the point $P_i$. Let
$\varphi$ contracts a curve $E$ different from $E_i$ for all $i$.
Then the pair $(S_{\min}\ni \overline P,\overline\Delta)$ is canonical, where
$\overline P=\varphi(E)$. It is easy to prove that
$(S_{\min}\ni \overline P,\overline\Delta)$ is analytically isomorphic either
$(\CC^2\ni 0,\frac13\{x^3+y^l=0\})$, where $l=3,4$, or $(\CC^2\ni
0,\frac13\{x^3+xy^3=0\})$, or $(\CC^2\ni
0,\frac13\{x^2y+y^4=0\})$.
\begin{lemma}
There exists a surface $S'$ such that
$\psi(E_i)$ is $(+1)$ non-singular rational curve for some $i$, where $\varphi\colon
S\stackrel{\psi}{\longrightarrow} S' \longrightarrow S_{\min}$.
\begin{proof}
If $\overline\Delta=\frac13\overline E_j$ then $p_a(\overline E_j)=28$.
The curve $\overline E_j$ must have two singular points of multiplicity
4 and eight singular points of multiplicity 3, we get a contradiction. Therefore there
exists a curve
$\overline E_j$ such that
$(\OO_{\PP^2}(5)-\overline E_j)$ is a nef divisor.
Take a resolution of the curve
$\overline E_j$ singularities. We obtain a curve with a self-intersection index
$\ge +1$. Q.E.D.
\end{proof}
\end{lemma}
We can assume that the linear system
$|\psi(E_i)|$ gives a birational morphism
$S'\to S_{\min}$ \cite[proposition
1.10]{KudLd} and $\overline\Delta=\frac13\overline
E_1+\frac13\overline E_2+\frac13\overline E_3$.
\par

Sorting out all variants of $\overline E_1$, $\overline E_2$, $\overline E_3$
on $S_{\min}\cong \PP^2$ the reader will easily prove that
there are four (--1) curves on $\widetilde S$ such that every curve
intersects all
$E_i$ and they are mutually disjoint.
Let us contract them $S\to S'$. We obtain a surface
$S'$ from case
3) of example \ref{ex1}.

\end{proof}
\end{theorem}

\begin{subs} The classification of log del Pezzo surfaces is very important
to study the three-dimensional extremal contractions and singularities,
because there is an induction from a
(local) three-dimensional contraction to
a two-dimensional log variety \cite{PrSh}, \cite{Sh2},
\cite{PrBound}.
Let us remark that in the result of induction
we obtain the log surfaces
$(S,D)$ such that
the divisor $-(K_S+D)$ is nef, big and the coefficients
of $D$ are standard. In order to get an effective classification,
the exceptional log del Pezzo surfaces
with $\delta=1,2$ are considered in more wide set of coefficients
-- $\M$ \cite{Sh2}, \cite{KudLd}. In this case the big condition is replaced, for
instance, on the requirement of existence of
$\QQ$-complement of $K_S+D$. The later allows to give the classification
of log Enriques surfaces with
$\delta=1,2$
\cite{KudLe1}, \cite{KudLe2}. Therefore, in the case $\delta=0$
the set $\SM$ will be extended to $\Phi_i$.
Now the main goal is to construct the models of $(S,D)$
with Picard number 1 or 2.
\end{subs}

\begin{definition} Let $(S,D)$ be an exceptional log del Pezzo surface with
$\delta(S,D)=0$, where $S$ is a rational surface. Then the pair
$(S,D)$ of type $\Phi_i$ is called {\it an exceptional log del Pezzo surface of
type} $\Phi_i$, where $i=2,3,4,5,6$.
\end{definition}

\begin{subs} Let $i=3,4,5,6$. Put
$$
\widehat D=\sum_{d_k\ge \frac{i-1}i}c(D_k)D_k+\sum_{d_k <
\frac{i-1}i}d_kD_k,
$$
where $c(D_k)=c(S,D-d_kD_k;D_k)$ is a log canonical threshold of a divisor
$D_k$ for the pair $(S,D-d_kD_k)$. By proposition \ref{two3}
the divisor $K_S+\widehat D$ is log canonical. By corollary
\ref{cor1} the divisor $-(K_S+\widehat D)$ is not nef.
Assume that $\rho(S)\ge 3$. Then there exists an exceptional curve
$E$ with $(K_S+\widehat D)\cdot E>0$ (see the proof of
theorem 4.1 from the paper \cite{Sh2}). Let
$\varphi\colon S\to S'$ be a contraction of $E$.
In contrast to the case
$\delta(S,D)\ge 1$ the birational morphism $\varphi$ can contract the curve
from $D$ with a coefficient $\ge \frac{i-1}i$.
\end{subs}

\begin{lemma}\label{two4}
Let $D=\alpha E+D^{\circ}$, $\alpha\ge \frac{i-1}i$ and
$P=\varphi(E)$. Then we have one of the following cases.
\begin{enumerate}
\item $(S'\ni P,\varphi(D))\an(\CC^2\ni
0,\frac12\{x=0\}+(\frac34+\varepsilon) \{x+y^3=0\})$, where
$\varepsilon<\frac1{60}$ and $i=4$. A morphism $\varphi$ is a weighted blow-up with
weights
$(3,1)$.
\item $(S'\ni P,\varphi(D))\an(\CC^2\ni 0,(\frac23+\varepsilon_1)\{x=0\}+
(\frac23+\varepsilon_2) \{x+y^2=0\})$, where
$\varepsilon_1+\varepsilon_2<\frac1{24}$ and $i=3$. A morphism
$\varphi$ is a weighted blow-up with
weights $(2,1)$.
\item
$(S'\ni P,\varphi(D))\an(\CC^2\ni
0,\frac12\{x=0\}+(\frac23+\varepsilon) \{x+y^l=0\})$, where $l=3,4$,
$i=3$ and $\varepsilon<\frac3{4l}-\frac16$. A morphism $\varphi$
is a weighted blow-up with
weights $(l,1)$.
\item $(S'\ni
P,\varphi(D))\an(\CC^2\ni
0,(\frac{i-1}i+\varepsilon)\{x^2+y^3=0\})$, where $i=5,4,3$,
$\varepsilon<\frac1{180}$ if $i=5$, $\varepsilon<\frac1{20}$ if
$i=4$, $\varepsilon<\frac1{12}$ if $i=3$. A morphism $\varphi$
is a weighted blow-up with
weights $(3,2)$.
\end{enumerate}
\begin{proof}
Note that $E$ must intersect a curve from
$D^{\circ}$ with a coefficient
$\ge \frac{i-1}i$. Let us consider the case when
$i=3$ and there exists a point on
$E$ with
$k=3$ from point 2) of proposition \ref{two3}. Then
\begin{gather*}
(K_S+D)\cdot E>(K_S+E+D^{\circ})\cdot
E=-2+\deg\Diff_E(0)+D^{\circ}\cdot E \ge \\
\ge -2+\frac32+\frac23>0.
\end{gather*}

A contradiction. Therefore $E$ has a coefficient 1 in $\widehat D$.
There are three possibilities.
\\
1). Assume that $\Diff_E(0)=0$. Then $(S'\ni P)\an (\CC^2\ni
0)/\ZZ_n(1,1)$. Since $(K_S+\widehat D)\cdot E>0$ and
$(K_S+D)\cdot E\le 0$ then simple calculations
show that the possibility is not realized.
\\
2). Assume that $\Diff_E(0)=\frac{k-1}kP$, where $k\ge 2$.
For example, consider the case
$i=3$. Since
$\deg\Diff_E(D^{\circ})<2$ and $\deg\Diff_E(\widehat D-E)>2$ then
$k=2,3,4,5$ and $P$ is a non-singular point of $S'$.
Moreover we have
\begin{enumerate}
\item[1)] $\Diff_E(D^{\circ})
=\frac12P+(\frac23+\varepsilon_1)P_1+(\frac23+\varepsilon_2)P_2$
or
\item[2)] $\Diff_E(D^{\circ})=
\frac{k-1}kP+\frac12P_1+(\frac23+\varepsilon)P_2$, where $k=3,4,5$.
\end{enumerate}
It remains to check that $-a(E,\varphi(D))<\frac{i}{i+1}$. We obtain cases
(2) and (3).
\\
3). Assume that
$\Diff_E(0)=\frac{k_1-1}{k_1}P_1+\frac{k_2-1}{k_2}P_2$, where
$k_1,k_2\ge 2$. Since $\deg\Diff_E(D^{\circ})<2$ then it can be assumed that
$k_1=2$. By direct calculations we obtain case (4).
\end{proof}
\end{lemma}

\begin{subs}\label{s1}
If $i=4,5,6$ then we repeat the procedure for $S'$ described above.
If $i=4$ then the case, when there is a point on an exceptional curve from case
(1) of lemma \ref{two4}, is impossible by the same argument as case 2)
of proposition \ref{two3} with $k=3$ (see the proof of lemma
\ref{two4}). As a result we get a surface $\overline S$ with
$\rho(\overline S)=1$, or a surface $\overline S$ with
$\rho(\overline S)=2$ and with structure of generically $\PP^1$-fibration.
\par
Now let $i=3$. Put $D'=\sum d_k'D_k'=\varphi(D)$.
Let us repeat the procedure described above.
If there is no a point on $D'_k$ from
case (2) of lemma \ref{two4},
the divisor $\widehat D'$ is defined as the divisor
$\widehat D$. Otherwise, put
$c(D'_k)=3/4$.
\par
Let $\varphi'\colon S'\to S''$ be a contraction of a curve $E'$ from
$D'$ with a coefficient $\ge 2/3$.
Two new cases can be appeared.
\\
1). The case, when there is a point on $E'$ from case (3) of lemma
\ref{two4}, is similarly impossible.\\
2). Assume that there is a point on $E'$ from case (2) of lemma
\ref{two4}. Then $\Diff_{E'}(\widehat D'-\frac34E')\ne \frac32Q$,
where $Q$ is a point of $E'$. Indeed, otherwise we have
$$
\left\{
\begin{array}{l}
(K_{S'}+\widehat D')\cdot E'=(K_{S'}+\frac34E'+3/4D'_2)\cdot E'=
-\frac12-\frac14{E'}^2>0\\
(K_{S'}+ D')\cdot
E'=(K_{S'}+(\frac23+\varepsilon_1)E'+(\frac23+\varepsilon_2)D'_2)
\cdot E'=\\ \hspace{5.6cm}
=-\frac23+2\varepsilon_2-(\frac13-\varepsilon_1){E'}^2\le 0.\\
\end{array}
\right.
$$
Since $\varepsilon_1+\varepsilon_2<\frac1{24}$ then this system of inequalities
is contradictorily. The simple calculations show that
the following new case is possible only.
\begin{equation}\tag{I}
(S''\ni P',\varphi'(D'))\an (\CC^2\ni
0,(\frac23+\varepsilon)\{x^2+y^5=0\}),
\end{equation}
where $\varepsilon<\frac1{120}$, $P'=\varphi'(E')$ and $\varphi'$ is a
weighted blow-up with weights
(2,1). Repeating the procedure we obtain the surface
$\overline S$ described above.
\end{subs}

\begin{subs}
Let $i=2$. Now the main problem is how much to increase the divisor
$D=\sum d_lD_l$ up to  $\widehat D=\sum
c_lD_l$, where $c_l\ge d_l$ for all $l$.
After this increase the divisor
$K_S+\widehat D$ must be log canonical, but not kawamata log terminal.
As before
there is an exceptional curve
$E$ with $(K_S+D)\cdot E>0$ under the condition
$\rho(S)\ge 3$. The corresponding morphism is denoted by
$\varphi\colon S\to S'$.
\par
Let us describe the construction of $\widehat D$. Let $D_t$ be a non-singular curve.
If there is a point on
$D_t$ from point 1) or 2) of proposition
\ref{two3} with $k\ge 4$ then $E\ne D_t$. Indeed, otherwise we have
$$
(K_S+D)\cdot D_t>(K_S+D_t+\sum_{l\ne t}d_lD_l)\cdot D_t\ge -2+
(\sum_{l\ne t}d_lD_l)\cdot D_t\ge 0.
$$
Therefore it is not important how much to increase the coefficient
$d_t$ in this case.
\par
Consider the remaining cases. If there is a point on
$D_t$ from point
1) or 2) of proposition \ref{two3} with $k=3$ then we put
$c_t=\frac23$.
\par
Consider the remaining cases.
If there is a point on
$D_t$ from point
1) or 2) of proposition \ref{two3} with $k=2$ then we put
$c_t=\frac34$.
\par
In the remaining cases we arbitrarily increase the other coefficients
$d_l$ up to maximal possible values.
\end{subs}

\begin{lemma}\label{two5}
Let $D=\alpha E+D^{\circ}$, $\alpha\ge \frac12$ and
$P=\varphi(E)$. Then we have one of the following cases.
\begin{enumerate}
\item $(S'\ni P,\varphi(D))\an(\CC^2\ni 0,(\frac12+\varepsilon)
\{x^3+y^4=0\})$, where $\varepsilon<\frac1{24}$. A morphism $\varphi$
is a weighted blow-up with weights $(1,1)$.
\item $(S'\ni
P,\varphi(D))\an(\CC^2\ni 0,(\frac12+\varepsilon_1)\{y=0\}+
(\frac12+\varepsilon_2) \{x^2+y^3=0\})$, where
$3\varepsilon_1+\varepsilon_2<\frac14$.
A morphism $\varphi$
is a weighted blow-up with weights $(1,1)$.
\item $(S'\ni
P,\varphi(D))\an(\CC^2\ni 0,(\frac12+\varepsilon)
\{x^2+y^{2k+1}=0\})$, where $\varepsilon<\frac3{8k+4}$.
 A morphism $\varphi$
is a weighted blow-up with weights
$(k,1)$.
\item
$(S'\ni P,\varphi(D))\an(\CC^2\ni
0,(\frac12+\varepsilon_1)\{x=0\}+
(\frac12+\varepsilon_2)\{y=0\}+(\frac12+\varepsilon_3)
\{x+y=0\})$, where
$\varepsilon_1+\varepsilon_2+\varepsilon_3<\frac16$.
 A morphism $\varphi$
is a weighted blow-up with weights
$(1,1)$.
\item
$(S'\ni P,\varphi(D))\an(\CC^2\ni
0,(\frac12+\varepsilon_1)\{x=0\}+ (\frac12+\varepsilon_2)
\{x+y^k=0\})$, where $\varepsilon_1+\varepsilon_2<\frac2{3k}$.
 A morphism $\varphi$
is a weighted blow-up with weights
$(k,1)$.
\end{enumerate}
\begin{proof}
Since the divisor $K_S+D$ is $\frac12$-log terminal then
the singularities of $S$ lying on $E$ are Du Val singularities
of type $\AAA_n$ (see lemma \ref{two3}).
For the same reason there is at most one singular point
on a curve
$E$.
\par
Assume that $E$ has the point of tangency of multiplicity 3 with a curve
$D_1$ from $D^{\circ}$. Since $(K_S+D)\cdot E\le 0$ then
$\Diff_E(D^{\circ}-d_1D_1)=0$ and we obtain case (1).
\par
If $E$ has the point of tangency of multiplicity 3 with a curve
$D_1$ from
$D^{\circ}$ then arguing as above in lemma
\ref{two4} we obtain cases
(2) and (3).
\par
For the remaining possibility we have cases (4) and (5).
\end{proof}
\end{lemma}

\begin{subs}\label{s2}
Let us repeat the procedure described above.
As a result of multiple procedure repetition two new
singularities can appear similarly to point
\ref{s1}:
\begin{equation}
\tag{II} (\CC^2\ni
0,(\frac12+\varepsilon_1)\{y=0\}+(\frac12+\varepsilon_2)\{x=0\}+
(\frac12+\varepsilon_3)\{x+y^k=0\}),
\end{equation}
where $\varepsilon_1+k(\varepsilon_2+\varepsilon_3)<\frac16$, and
\begin{equation}
\tag{III} (\CC^2\ni 0,(\frac12+\varepsilon_1)\{y=0\}+
(\frac12+\varepsilon_2)\{x^2+y^{2k+1}=0\}),
\end{equation}
where $\varepsilon_1+2k\varepsilon_2<\frac16$.
\end{subs}

The results above-mentioned allow to define the model
of exceptional log del Pezzo surface of type $\Phi_i$ (cf.
\cite[\S 5]{Sh2}).
\begin{definition}\label{maindef}
Let $(S,D)$ be an exceptional log del Pezzo surface of type
$\Phi_i$ $(i=2,3,4,5,6)$ except the following points (if
$i=6$ then there are no exceptions).
\begin{enumerate}
\item If $i=5$ then see point (4) of lemma \ref{two4}.
\item If
$i=4$ then see points (1), (4) of lemma \ref{two4}.
\item If $i=3$ then see points
(2), (3), (4) of lemma \ref{two4} and case
(I) of point \ref{s1}.
\item If $i=2$ then see points (1), (3),
(5) of lemma \ref{two5}. It is possible with another restrictions on the values of
$\varepsilon$, $\varepsilon_1$, $\varepsilon_2$. Also see cases
(II) and (III) of point \ref{s2}.
\end{enumerate}
Then the pair $(S,D)$ is called {\it a model of type $\Phi_i$} if
one of the following two conditions is satisfied.
\begin{enumerate}
\item[A).] $\rho(S)=1$.
\item[B).] $\rho(S)=2$, the cone $\NE(S)$ is generated by two extremal rays
$R_1$ and $R_2$. The ray $R_1$ gives generically
$\PP^1$-fibration. If the ray $R_2$ gives a birational contraction
of a curve $E$ then $E$ is a component of divisor $D$ with a coefficient
$\ge \frac{i-1}i$.
\end{enumerate}
\end{definition}

\begin{remark} Let us remark that in the model of type
$\Phi_i$ definition the condition, that the divisor $K_S+D$ is
$\frac1i$-log terminal, is not fulfilled at the non-singular points of
surface $S$ only.
\end{remark}

\begin{subs}
The very important problem is to classify the models of type
$\Phi_i$. When the model classification of type
$\Phi_2$ is finished it is remained to describe the exceptional surfaces
$S$ such that the divisor $K_S$ is
$\frac12$-log terminal
and there exists
a 2-complement $\Theta=\sum\frac12\Theta_i$ of $K_S$ (see theorem
\ref{main3}). This completes the classification of exceptional log
del Pezzo surfaces
(see \cite{KudLd}, \cite{Sh2}) and allows to describe
log Enriques surfaces completely
(see
\cite{KudLe1}, \cite{KudLe2}).
\end{subs}


\begin{thebibliography}{99}
\bibitem{Bl}
\textit{Blache R.} The structure of l.c. surfaces of Kodaira
dimension zero // J. Algebraic Geom. 1995. V. 4. P. 137--179.

\bibitem{KeM}
\emph{Keel S., McKernan J.} Rational curves on quasi-projective
surfaces // Memoirs AMS 1999. V. 140. no. 669.

\bibitem{Koetal}
\emph{Kollar J. et al} Flips and abundance for algebraic
threefolds // Ast\'erisque 1992. V. 211.

\bibitem{Kud2}
\emph{Kudryavtsev S.~A.} Pure log terminal blow-ups // Math.
Notes. 2001. V. 69. No. 6. P. 814--819.

\bibitem{Kud3}
\emph{Kudryavtsev S.~A.} Classification of three-dimensional
exceptional log canonical hypersurface singularities.
{\small I} // Russian Acad. Sci. Izv. Math.
2002. V. 66. No. 5. P. 949-1034.

\bibitem{Kud4}
\emph{Kudryavtsev S.~A.} Classification of three-dimensional
exceptional log canonical hypersurface singularities.
{\small II} // e-print math.AG/0202105.

\bibitem{KudLe2}
\emph{Kudryavtsev S.~A.} Classification of logarithmic Enriques surfaces with
$\delta=2$ // Math. Notes 2002. V. 72 no. 5.  P. 660-666.

\bibitem{KudLe1}
\emph{Kudryavtsev S.~A.} Classification of logarithmic Enriques surfaces with
$\delta=1$ // appear in Math. Notes 2003; e-print math.AG/0207270

\bibitem{KudLd}
\emph{Kudryavtsev S.~A.} Classification of exceptional log del Pezzo surfaces with
$\delta=1$ // Russian Acad. Sci. Izv. Math.
2003. V. 67. No. 3; e-print math.AG/0207269.

\bibitem{PrRed}
\emph{Prokhorov Yu.~G.} Mori conic bundles with a reduced log
terminal boundary// J. Math. Sci. 1999. V. 91. no. 1. P.
1051--1059.

\bibitem{PrBound}
\emph{Prokhorov Yu.~G.} Boundedness of non-birational extremal
contractions // Internat. J. Math. 2000. V. 11. P. 393--411.

\bibitem{PrSh}
\emph{Prokhorov Yu.~G., Shokurov V.~V.}
The first main theorem on complements:
from global to local // Russian Acad. Sci. Izv. Math.
2001. V. 65. No. 6. P. 1169-1196.

\bibitem{PrLect}
\emph{Prokhorov Yu.~G.} Lectures on complements on log surfaces //
MSJ Memoirs V. 10. 2001.

\bibitem{Hart}
\emph{Hartshorne R.} Algebraic Geometry // Springer  1977.

\bibitem{Sh1}
\emph{Shokurov V.~V.} $3$-fold log flips // Russian Acad. Sci.
Izv. Math. 1993. V. 40. P. 93--202.

\bibitem{Sh2}
\emph{Shokurov V.~V.} Complements on surfaces // J. of Math. Sci.
2000. V. 102. no. 2. P. 3876--3932.

\bibitem{Zh1}
\textit{Zhang D.-Q.} Logarithmic Enriques surfaces //J. Math.
Kyoto Univ. 1991. V. 31. P. 419--466.

\end{thebibliography}
\end{document}